\documentclass[11pt]{amsart}

\usepackage{amsmath}
\usepackage{xypic}
\usepackage{amssymb}
\usepackage{hyperref}

\newtheorem{theorem}{Theorem}[section]
\newtheorem{proposition}[theorem]{Proposition}
\newtheorem{lemma}[theorem]{Lemma}
\newtheorem{definition}[theorem]{Definition}

\newtheorem{corollary}[theorem]{Corollary}
\newtheorem{remark}[theorem]{Remark}

\begin{document}

\def\be{\begin{equation}}
\def\ee{\end{equation}}

\def\ra#1{\mathop{\longrightarrow}\limits^{#1}}
\def\rab#1{\mathop{\longrightarrow}\limits_{#1}}
\def\lab#1{\mathop{\longleftarrow}\limits_{#1}}
\def\la#1{\mathop{\longleftarrow}\limits^{#1}}
\def\da#1{\downarrow\rlap{$\vcenter{\hbox{$\scriptstyle#1$}}$}}
\def\ua#1{\uparrow\rlap{$\vcenter{\hbox{$\scriptstyle#1$}}$}}
\def\sea#1{\mathop{\searrow}\limits^{#1}}
\def\nea#1{\mathop{\nearrow}\limits^{#1}}
\def\swa#1{\mathop{\swarrow}\limits^{#1}}
\def\dirlim#1{\lim_{\textstyle{\mathop{\rightarrow}\limits_{#1}}}}
\def\invlim#1{\lim_{\textstyle{\mathop{\leftarrow}\limits_{#1}}}}

\def\AA{\mathbb{A}}
\def\DD{\mathbb{D}}
\def\QQ{\mathbb{Q}}
\def\ZZ{\mathbb{Z}}
\def\NN{\mathbb{N}}
\def\GG{\mathbb{G}}
\def\WW{\mathbb{W}}
\def\RR{\mathbb{R}}
\def\CC{\mathbb{C}}
\def\FF{\mathbb{F}}
\def\SS{\mathbb{S}}
\def\bzero{\mathbf{0}}
\def\bone{\mathbf{1}}

\def\bAlg{\mathbf{Alg}}
\def\bCM{\mathbf{ComMon}}
\def\bCA{\mathbf{ComAlg}}
\def\bSet{\mathbf{Set}}
\def\bAb{\mathbf{Ab}}
\def\bMod{\mathbf{Mod}}
\def\bC{\mathbf{C}}
\def\bA{\mathbf{A}}
\def\bD{\mathbf{D}}
\def\bL{\mathbf{L}}
\def\bRMod{\mathbf{RMod}}
\def\bLMod{\mathbf{LMod}}

\def\cP{\mathcal{P}}

\def\qed{$\,\square$}

\def\tensor{\otimes}
\def\dsum{\oplus}
\def\Tensor{\mathop\bigotimes}
\def\Dsum{\mathop\bigoplus}

\def\bR{\mathbf{R}}

\def\cW{\mathcal{W}}
\def\cC{\mathcal{C}}
\def\cF{\mathcal{F}}
\def\cG{\mathcal{G}}
\def\cO{\mathcal{O}}
\def\cL{\mathcal{L}}
\def\cM{\mathcal{M}}
\def\cN{\mathcal{N}}

\def\Hom{\mathrm{Hom}}
\def\c{\mathrm{c}}
\def\inj{\mathrm{in}}
\def\sgn{\mathrm{sgn}}
\def\ch{\mathrm{ch}\,}
\def\Ab{\mathrm{Ab}}
\def\id{\mathrm{id}}
\def\Der{\mathrm{Der}}
\def\pr{\mathrm{pr}}
\def\Tor{\mathrm{Tor}}
\def\Ext{\mathrm{Ext}}
\def\Sym{\mathrm{Sym}}

\def\intersect{\mathop\bigcap}
\def\iso{\cong}
\def\cotensor{\square}
\def\tc{\textstyle{\coprod}}

\def\Rb{{R^\bullet}}
\def\Sb{{S^\bullet}}
\def\ssin{\hbox{\scriptsize in}}
\def\hG{\widehat{G_0}}

\def\oC{\overline{C}}
\def\oI{\overline{I}}

\def\wQQ{\widetilde{\mathbb{Q}}}
\def\wZZ{\widetilde{\mathbb{Z}}}
\def\wK{\widetilde{K}}

\title
{The Andr\'e-Quillen cohomology of commutative monoids}

\author{Bhavya Agrawalla}
\address{Department of Mathematics, Massachusetts Institute of Technology,
Cambridge, MA 02139}
\email{bhavya@mit.edu}

\author{Nasief Khlaif}
\address{Department of Mathematics, Birzeit University,
Ramallah, Palestine}
\email{nkhlaif@birzeit.edu}

\author{Haynes Miller}
\address{Department of Mathematics, Massachusetts Institute of Technology,
Cambridge, MA 02139}
\email{hrm@math.mit.edu}



\begin{abstract}
We observe that Beck modules for a commutative monoid are exactly
modules over a graded commutative ring associated to the monoid.
Under this identification,
the Quillen cohomology of commutative monoids is a special
case of Andr\'e-Quillen cohomology for graded commutative rings,
generalizing a result of Kurdiani and Pirashvili.
To verify this we develop the necessary grading formalism.
The partial cochain complex developed by Pierre Grillet for computing
Quillen cohomology appears as the start of a modification
of the Harrison cochain complex suggested by Michael Barr.
\end{abstract}


\setcounter{equation}{0}

\maketitle

In his book {\em Homotopical Algebra} \cite{quillen},
Daniel Quillen described a homotopy theory of simplicial objects in any of
a wide class of universal algebras, and  corresponding theories of homology
and cohomology. Quillen homology is defined as derived functors of an
abelianization functor, and in many cases can be computed using a cotriple
resolution \cite{barr-beck}. 
Coefficients for these theories are ``Beck modules,'' that is,
abelian objects in a slice category. The case of commutative rings was 
studied at the same time by Michel Andr\'e \cite{andre}.

An example of such an algebraic theory, one of long standing and 
increasing importance, is provided by the category $\bCM$ of commutative 
monoids. The prime exponent of the study of commutative monoids
has for years been Pierre Grillet 
\cite{grillet-74,grillet-95,grillet-97,grillet-01,grillet-22,grillet-22b}
(but see also \cite{calvo} and \cite{kurdiani-pirashvili} for example).
Among other things, Grillet
provided the beginning of a small cochain complex, based on multilinear
maps subject to certain symmetry conditions, whose cohomology he showed to be
isomorphic in low dimensions to the Quillen cohomology $HQ^*_{CM}(X;M)$ 
of the commutative monoid $X$ with coefficients in a Beck module $M$ for $X$;
and in \cite{grillet-21} a corresponding resolution in Beck modules
was developed. These results are surprising, since Quillen cohomology
is defined by means of a simplicial resolution and does not generally
admit such an efficient computation.

It is well-known (\cite{leech,grillet-74} and \cite[p.~29]{wells}) and easy
to see that the category of (left) Beck modules for $X$, $\bLMod_X$,
is equivalent to the category of 
covariant functors from the ``Leech category'' $L_X$ to the category
$\bAb$ of abelian groups. The Leech category has object set $X$; 
$L_X(x,y)=\{z:y=z+x\}$ (writing the monoid structure additively); 
and composition is given by addition in the commutative monoid. 

In this paper, we observe that Grillet's construction is in fact
subsumed by the theory of Harrison cohomology of commutative rings,
slightly augmented as suggested by Michael Barr,
once this theory has been extended to the graded context. As pointed
out by Bourbaki \cite[Ch. 2 \S11]{bourbaki-algebra-1}, 
one can speak of rings graded
by a commutative monoid: an $X$-graded object in a category $\bC$ is an 
assignment of an object $C_x\in\bC$ for each $x\in X$. If $\bC$ is the
category $\bMod_K$ of modules over some commutative ring $K$, there is
a natural symmetric monoidal structure on the category $\bMod_K^X$ 
of $X$-graded $K$-modules, and we may define $X$-graded $K$-algebras,
and modules over them, accordingly. 

The first observation, simple enough as to need no proof, is that there is
a natural $X$-graded commutative $K$-algebra $\wK X$ in which, 
for each $x\in X$, $(\wK X)_x$ is the free $K$-module generated by an 
element we will write $1_x$, with the evident
unit and multiplication. This is the ``$X$-graded monoid $K$-algebra'' of $X$. 

The next observation, equally simple, is that the category of Beck $X$-modules
is equivalent to the category of ($X$-graded) left modules over $\wZZ X$. 

These two observations bring into play the entire highly developed homological
theory of commutative rings. Our first main result \ref{cm-ca} is that 
\[
HQ_{CM}^*(X;M)=HQ_{CA}^*(\wZZ X;M)
\]
where the right hand term denotes the well-studied Andr\'e-Quillen cohomology
\cite{quillen-com-alg,andre,weibel}, 
extended to the graded context. This generalizes an observation of
Kurdiani and Pirashvili \cite{kurdiani-pirashvili}, who considered the case 
of Beck modules pulled back from the trivial monoid, in which case one 
arrives at the Andr\'e-Quillen cohomology of $\ZZ X$ as an ungraded
commutative ring.

Andr\'e-Quillen cohomology is of course hard to compute, but there are 
well known approximations to it. One such approximation is given by
Harrison cohomology \cite{harrison,gerstenhaber-schack,barr,whitehouse}
$Harr^*(A;M)$. 
This theory is most neatly expressed by restricting to Hochschild cochains
that annihilate shuffle decomposables; or, equivalently, to cochains that
satisfy appropriate ``partition'' symmetry conditions.  
This characterization was apparently
suggested by Mac Lane, and was adopted in \cite{harrison}, but Harrison's 
original invariance property involved a different characterization
of the same symmetry conditions, using ``monotone'' permutations. 
The equivalence of these two definitions can be found as Corollary 4.2 in
\cite{gerstenhaber-schack}.

This approximation definitely breaks down in finite characteristic:
The Andr\'e-Quillen
cohomology of a polynomial algebra vanishes in positive dimensions,
but Michael Barr showed \cite{barr} 
that the Harrison cohomology of the polynomial algebra over a field of
characteristic $p$ is nonzero in dimension $2p$.
Barr himself proposed a variant of the Harrison construction, restricting
Hochschild cochains that not only by vanish on shuffle decomposables but
also on divided powers. This overcomes the obstacle in dimension $2p$, but
in her thesis \cite{whitehouse} Sarah Whitehouse proved
that this variant also fails to give Andr\'e-Quillen
cohomology, by showing that Barr cohomology in dimension 5 does not vanish
on $\FF_2[x]$.

Our second observation is that 
exactly the same monotone symmetry and divided-power annihilation
conditions occur in the partial complex 
described by Grillet; Grillet's partial complex is precisely
the beginning of the ``Barr complex'' for the graded monoid
algebra. In later work \cite{grillet-22}, Grillet proposed that this complex
correctly computes Quillen cohomology in higher dimensions as well, but
Whitehouse's counterexample shows that this conjecture
fails at least in dimension 5.

Grillet's identification of his cohomology with the Quillen cohomology
of commutative monoids goes well beyond what seems to be known about the
relationship between Harrison cohomology and Andr\'e-Quillen cohomology
in general, and suggests a variety of questions.

We note that the observation that Beck modules over a 
commutative monoid are just graded modules over its graded monoid algebra
suggests that the description of Quillen homology for commutative monoids
carried out in \cite{kurdiani-pirashvili} is in fact a special case of a
graded extension of Pirashvili's earlier work \cite{pirashvili}.

We begin in \S1 with a recollection of Quillen homology, along with the
cotriple resolution that may be used to compute it. \S2 sets out some 
elementary facts about $X$-gradings, and in 
\S3 we explain how the grading behaves in homological algebra. Change
of grading is explained in \S4. The next section is the core of the work,
proving the identification of Quillen homology and cohomology 
for commutative monoids with that of certain $X$-graded commutative rings. 
We then turn to interpreting
the work of Grillet. This requires developing the Hochschild complex  
with its shuffle product and divided power structure (and we provide
some new general information about the latter), and the various indecomposable 
quotients occuring in the definitions of Harrison and Barr homology. 
In \S9, we review the motivating work of Pierre Grillet and relate it
to Harrison and Barr cohomology.

\smallskip
\noindent
{\bf Acknowledgements.}
We are grateful to Pierre Grillet for forwarding us, in response to a letter
from us outlining the results presented here, an early copy of
a paper in which a similar story is worked out. He uses somewhat different
language -- his ``multi'' objects are our graded objects -- but he did not
make the connection with Harrison homology that we establish here.

This work was carried out under the auspices of a program, supported by
MIT's Jameel World Education Laboratory, designed to foster
collaborative research projects involving students from MIT and Palestinian
universities.
We acknowledge with thanks the contributions made by early participants
in this program -- Mohammad Damaj and Ali Tahboub of Birzeit University 
and Hadeel AbuTabeekh of An-Najah National University -- as well as the 
support of Palestinian faculty -- Reema Sbeih and Mohammad Saleh at Birzeit
and Khalid Adarbeh and Muath Karaki at NNU. We thank Professor Victor Kac
for pointing out to us the relevance of \cite{gerstenhaber-schack}.
The first author acknowledges support by the MIT UROP office. 

Finally, we thank the referee for such a careful reading of the document.

\section{Quillen homology and Quillen cohomology}

In \cite{quillen-com-alg}, Daniel Quillen proposed a uniform definition of the
``homology'' and ``cohomology'' of objects in a very general class of
categories. The marquis example was that of commutative algebras, but
his definition applies much more generally and subsumes many of the ad hoc
definitions that were already in use at the time and have been considered
subsequently.

Quillen proposed that the construction of ``homology'' should be a special
case of a general procedure for deriving a functor. One of the motivations
for his development of the theory of ``model categories''
\cite{quillen,hovey,hirschhorn,goerss-jardine}
was precisely to provide a context for
defining derived functors of non-additive functors. This theory
``internalizes'' homological algebra, in the sense that objects playing
the role of projective resolutions (called ``cofibrant objects'')
and maps playing the role of quasi-isomorphisms (called ``weak equivalences'')
exist in the category (rather than in some auxiliary
category such as a category of chain complexes).
One of the axioms asserts that an object admits a weak equivalence from a
cofibrant object; this ``cofibrant replacement'' plays the role of a
projective resolution.

In \cite[II\S4]{quillen}, Quillen establishes the existence of a model
structure on the category of simplicial objects over any one of a very
general class of categories with suitable mild properties. We refer to
Quillen's book or \cite{frankland} for definitions.

\begin{theorem}{\cite[II\S4, Theorem 4]{quillen}}
Let $\bC$ be any cocomplete category admitting a set $\cP$ of small projective
generators (a ``quasi-algebraic category'' in the language of
\cite{frankland}).
Then there is a model structure on the
category of simplicial objects over $\bC$ in which the weak equivalences are
the morphisms $f$ such that $\bC(P,f)$ is a weak equivalence of simplicial
sets for all $P\in\cP$.
\end{theorem}

All normally occurring categories of universal algebras satisfy these 
assumptions.

This model structure on $s\bC$ allows us to define derived functors for any
functor $E:\bC\to\bA$, where $\bA$ is an abelian category:
For any $X$ in $\bC$,
let $P_\bullet\to X$ be a cofibrant replacement and define
\[
L_nE(X)=\pi_n(EP_\bullet)\,.
\]
See \cite{frankland} for an elaboration of the naturality of this construction.
The fact that
any cofibrant replacement can be used is contained, for example, in
\cite[Proposition 3.9]{frankland}.

Explicit cofibrant replacements can often be constructed as a ``cotriple
resolution'' \cite{barr-beck}. An adjoint pair
\[
\bC\rightleftarrows\bD
\]
defines a triple $F$ on $\bC$ and a cotriple $G$ on $\bD$; see
\cite{eilenberg-moore,beck}. For example the free commutative
monoid functor $\NN$ is left adjoint to the forgetful functor --
\[
\NN:\bSet\rightleftarrows\bCM:u
\]
-- and this adjoint pair defines triple $u\NN$ on $\bSet$ and a cotriple
$\NN u$ on $\bCM$. The commutative monoid of 
natural numbers (also denoted by $\NN$) is a small
projective generator for $\bCM$, and $u(S)=\bCM(\NN,S)$.

A cotriple $G$ on $\bC$ determines a functor $G_\bullet$ to the category
$s\bC$ 
of simplicial objects over $\bC$, augmented to the identity functor: this is
the ``cotriple resolution'': see \cite{barr-beck} or \cite[Chapter 8]{weibel}.
To relate it to the model category structure we need a further restriction:
A category $\bC$ is {\em algebraic} (in Frankland's sense) if it is
quasi-algebraic and Barr-exact 
(\cite[p. 35]{barr-book},\cite[p. 91]{frankland}).

\begin{proposition}{\cite[p.~ 69]{quillen-com-alg}}
    Let $\bC$ be an algebraic category and $\cP$ a set of small projective
    generators. 
    When the cotriple $G$ is associated to the adjoint pair
  $\bC\rightleftarrows\bSet^\cP$, the cotriple resolution $G_\bullet A\to A$
  serves as a cofibrant replacement of $A$ (regarded as a constant simplicial
  object) in the Quillen model structure on $s\bC$.
\end{proposition}

This allows one to calculate the derived functors for any functor $E:\bC\to\bA$
to an abelian category:
\[
L_nE(C)=\pi_n(E(G_\bullet C))=H_n(\ch E(G_\bullet C))
\]
where $\ch$ denote formation of the chain complex associated to a simplicial
object in an abelian category; see \cite[Theorem 8.4.1]{weibel}, for example. 

Quillen's definition \cite{quillen-com-alg} of homology and cohomology
of an object in a category $\bC$ involves the notion of Beck modules.

\begin{definition}{\cite[Definition 5]{beck}}
    A {\em Beck module} over an object $A$ in $\bC$ is
  an abelian object in the slice category $\bC/A$.
\end{definition}

With the evident morphisms, Beck $A$-modules form a category $\bLMod_A$.
The terminal object in $\bLMod_A$ is the identity map $\id_A:A\downarrow A$
with its unique abelian structure.
If $\bC$ is quasi-algebraic, so is $\Ab(\bC/A)$ \cite[3.40]{frankland}.
Under mild additional conditions $\Ab(\bC/A)$ is an abelian category:

\begin{proposition}{\cite[p. 69]{quillen-com-alg}
    \footnote{But Quillen inadvertently omits the exactness condition.},
  \cite[Chapter 2, Theorem 2.4]{barr-book}}
  Let $\bC$ be an algebraic category and $A\in\bC$.
  Then both $\bC/A$ and $\Ab(\bC/A)$ are algebraic; $\Ab(\bC/A)$ is
  abelian; and the forgetful functor $\Ab(\bC/A)\to\bC/A$ has a left
  adjoint $\Ab_A:\bC/A\to\Ab(\bC/A)$.
\end{proposition}

We are now in position to recall the following definition.

\begin{definition}{\cite{quillen}}
  Let $\bC$ be an algebraic category. 
The {\em Quillen homology} of an object $A$ in $\bC$ is the sequence
of Beck $A$-modules given by
\[
HQ_n(A)=L_n\Ab_A(\id_A)=\pi_n(\Ab_A(P_\bullet))
\]
where $P_\bullet\to A$ is a cofibrant replacement regarded as an object in
$s\bC/A$. 
The {\em Quillen cohomology} of $A$ with coefficients in a Beck $A$-module
is the sequence of abelian groups
\[
HQ^n(A;M)=H^n(\Hom_A(\ch P_\bullet,M))\,.
\]
\end{definition}
In terms of the cotriple resolution,
\begin{gather*}
  HQ_n(A)=\pi_n(\Ab_A(G_\bullet A))\\
  HQ^n(A;M)=H^n(\Hom_A(\Ab_A(G_\bullet A),M))\,.
\end{gather*}

For example \cite[\S4]{quillen-com-alg},
when $\bC$ is the category $\bCA_K$ of commutative $K$-algebras,
the category of Beck $A$-modules is equivalent to the category of left 
$A$-modules: An abelian object over $A$, $p:B\downarrow A$, first of all
has a section, the ``zero-section,'' which provides an identification of
$K$-modules $B\iso A\oplus M$ where $M$ is the kernel of $p$ as an $A$-module.
The abelian structure forces the multiplication on $A\oplus M$ to be given 
by $(a,m)(b,n)=(ab,an+bm)$. Under this identification, a section of 
$A\oplus M\downarrow A$ in $\bCA_K$ is given by $a\mapsto(a,da)$ where
$d\in\Der_K(A,M)$. The abelianization of $\id:A\downarrow A$ is the 
$A$-module such that $\Hom_A(\Ab_A(A),M)=\Der_K(A,M)$; that is, 
$\Ab_A(A)$ is the usual module $\Omega_{A/K}$ of K\"ahler differentials.
More generally, for $B\downarrow A$ in $\bCA_K/A$ 
\[
\Ab_A(B)=A\tensor_B\Omega_{B/K}\,.
\]
So in that case we have the ``cotangent complex''
\[
\bL_{A/K}=A\tensor_{G_\bullet A}\Omega_{G_\bullet A/K}\,.
\]
The Andr\'e-Quillen homology is its homotopy --
\[
HQ_n(A)=\pi_n(\bL_{A/K})
\]
-- and the Andr\'e-Quillen cohomology is
\[
HQ^n(A;M)=H^n(\Hom_A(\ch\bL_{A/K},M))\,.
\]
In this case we can also define homology with coefficients
in an $A$-module $M$:
\[
HQ_n(A;M)=\pi_n(\bL_{A/K}\tensor_AM)\,.
\]

\section{Gradings}
\label{sec-gradings}
We begin by setting up some structure on categories of objects graded
over a commutative monoid. 

Let $X$ be a commutative monoid, which we will write additively. 
Following Bourbaki \cite[Ch. 2 \S11]{bourbaki-algebra-1}, we say that 
an $X$-graded object $C_\bullet$ in a category $\bC$
is a choice of object $C_x$ of $\bC$ for each $x\in X$. 
Write $\bC^X$ for the category of $X$-graded objects in $\bC$.
A morphism $C_\bullet\to C'_\bullet$ is a morphism $C_x\to C'_x$ for
each $x\in X$. A functor $F:\bC\to\bD$ induces
$F^X:\bC^X\to\bD^X$, and an adjunction between $E$ and $F$ induces
a canonical adjunction between $E^X$ and $F^X$. 

Structure on $\bC$ often induces structure on $\bC^X$. For example suppose
that $(\bC,\bone,\tensor,c)$ is a closed symmetric
monoidal category \cite[\S6.1]{borceaux}.
Assuming that $\bC$ has coproducts of large enough sets of objects,
there is then a canonical symmetric monoidal structure on $\bC^X$, 
which is also closed
if $\bC$ has products of large enough sets of objects, in which 
\begin{gather*}
\bone_x=\begin{cases}
\bone & \text{for}\,\, x=0\\
\bzero & \text{for}\,\, x\neq0
\end{cases} \\
(C_\bullet\tensor D_\bullet)_z=\coprod_{x+y=z}C_x\tensor D_y\,.
\end{gather*}
The symmetry 
$c:(C_\bullet\tensor D_\bullet)_z\to(D_\bullet\tensor C_\bullet)_z$
is such that for all $x,y$ with $x+y=z$,
\[
c\circ\inj_{x,y}=\inj_{y,x}\circ c_{C_x,D_y}
\]
where $c_{C,D}:C\tensor D\to D\tensor C$ is the symmetry in $\bC$.

One may consider commutative monoids with respect to this
commutative monoid structure. 
For example, a commutative monoid in the symmetric monoidal category
$\bSet^X$ consists of a set $T_x$ for each $x\in X$ together with an
element $1\in T_0$ and maps $\mu:T_x\times T_y\to T_{x+y}$ satisfying
evident conditions. This is to be distinguished from an $X$-graded
commutative monoid, an object of $\bCM^X$! 

We can then describe the free commutative $\tensor$-monoid generated 
by an $X$-graded set $T$, denoted $\NN_X T$, 
by the following sequence of observations. 

\noindent{\bf(1)} 
The free commutative $\tensor$-monoid $\NN_X t$ with a single 
generator $t$ in degree $x$ has
\[
(\NN_X t)_y=\{k\in\NN:kx=y\}\,,
\]
and the product of $k\in(\NN_X t)_y$ with $k'\in(\NN_X t)_{y'}$ is 
$k+k'\in(\NN_X t)_{y+y'}$.

\noindent{\bf(2)} 
An $X$-graded set $T$ is {\em finite} if $\sum_x|T_x|<\infty$.
The free commutative $\tensor$-monoid generated by a finite $X$-graded set
$T$ is the tensor product of $\NN_X t$'s as $t$ runs over $T_x$, $x\in X$. 

\noindent{\bf(3)} 
Any $X$-graded set is the direct limit of a filtered family of 
finite $X$-graded sets, and the free commutative $\tensor$-monoid 
functor commutes with filtered colimits.
 
Let $K$ be a commutative ring. The category of 
$X$-graded $K$-modules $\bMod_K^X$
admits a symmetric monoidal structure given by the
``graded tensor product,'' with 
\[
(A_\bullet\tensor_K B_\bullet)_z=\bigoplus_{x+y=z}A_x\tensor_K B_y
\]
and unit given by the $X$-graded $K$-module with $K$ in degree $0$ and $0$ 
in all other degrees. The symmetry sends $x\tensor y$ to $y\tensor x$.
An $X$-graded $K$-algebra is a monoid for this tensor 
product. Once again, beware of this use of language; this is not an 
$X$-graded object in $\bCA_K$. Write $\bCA(\bMod_K^X)$ for this category. 
A Beck module for the commutative $X$-graded $K$-algebra $A_\bullet$ is an 
action of this monoid.

For $T\in\bSet^X$ the free commutative $X$-graded
$K$-algebra generated by $T\in\bSet^X$ is given in degree $x$
by the free $K$-module generated by $(\NN_X T)_x$. This
provides us with an adjoint pair
\[
F_X:\bSet^X\rightleftarrows\bCA(\bMod_K^X):u_X\,.
\]

The relationship with the Leech category (section \ref{sec-cm} below)
suggests that rather than defining a right $A_\bullet$-module
as a right action, we should say this:
\begin{definition} A right $A_\bullet$-module is an 
$X$-graded $K$-module $M^\bullet$ together with homomorphisms
\[
\varphi_{x,y}:M^{x+y}\tensor A_y\to M^x
\]
such that $\varphi_{x,0}(m\tensor1)=m$ and
\[\xymatrix{
M^{x+y+z}\tensor A_z\tensor A_y \ar[r]^{1\tensor\mu_{z,y}}
\ar[d]^{\varphi_{x+y,z}\tensor1} & 
M^{x+y+z}\tensor A_{z+y=y+z} \ar[d]^{\varphi_{x,y+z}} \\
M^{x+y}\tensor A_y \ar[r]^{\varphi_{x,y}} & M^x 
}\]
commutes. 
\end{definition}

Write $\bRMod_{A_\bullet}$ for the category of right $A_\bullet$-modules.
If $X$ has inverses, so is in fact an abelian group, the category
$\bRMod_{A_\bullet}$ is equivalent to the category of 
a right $A_\bullet$-modules in the usual sense, 
using ``lower indexing'' $M_x=M^{-x}$.

Let $N_\bullet$ be a left $A_\bullet$-module and $M^\bullet$ a right
$A_\bullet$ module. Their {\em tensor product} over $A_\bullet$,
$M^\bullet\tensor_{A_\bullet}N_\bullet$, is the $K$-module
defined as the coequalizer of the two maps
\[
f,g:P=\bigoplus_{x,y}M^{x+y}\tensor A_y\tensor N_x\rightrightarrows 
\bigoplus_z M^z\tensor N_z\,.
\]
Each of these maps is defined by giving the composite with an inclusion
$\inj_{x,y}:M^{x+y}\tensor A_y\tensor N_x\to P$:
\begin{align*}
f\circ\inj_{x,y}&=\inj_{x+y}\circ(1\tensor\varphi_{y,x})\,,\\
g\circ\inj_{x,y}&=\inj_x\circ(\varphi_{x,y}\tensor1)\,.
\end{align*}

\section{Graded homological algebra}

From now on we will write just $A$ rather than $A_\bullet$ and so on.
Let $A$ be an $X$-graded $K$-algebra. For each $x\in X$, there is 
a left $A$-module $P^x$ together with $\iota\in P^x_x$ such 
that for any left $A$-module $N$ the map
\[
\Hom_{A}(P^x,N)\to N_x\,,\quad f\mapsto f(\iota)
\]
is an isomorphism. Explicitly, 
\begin{equation}\label{pxy}
P^x_y=\bigoplus_{x+z=y}A_z
\end{equation}
and $\iota$ is the image of $1\in A_0$ under $\inj_0:A_0\to P^x_x$. Given
$n\in N_x$, the corresponding map $\widehat n:P^x\to N$ 
is defined by 
\[
\widehat n\circ\inj_z(a)=an\,,\quad a\in A_z\,.
\]

\begin{lemma} The set $\{P^x:x\in X\}$ is a generating set of small
  projective $A$-modules.
  \end{lemma}
\begin{proof}
For each $x\in X$, the $A$-module $P^x$ is projective since $N\mapsto N_x$ 
is an exact 
functor. An object is small if the functor it co-represents preserves
filtered colimits. For $P^x$ this is clear since colimits in $\bLMod_A$ 
are computed
component-wise. For any $A$-module $N$, the map
\[
\bigoplus_{x\in X}\bigoplus_{n\in N_x}P^x\to N\,,
\]
given by $\widehat n$ on the component indexed by $(x,n)$, is an epimorphism;
this shows that $\{P^x:x\in X\}$ is a generating set.
\end{proof}

It follows that any projective $A$-module is a retract of a direct sum of
$P^x$'s.

The account of Quillen homology and cohomology 
given above goes through in the graded
context without essential change. For an $X$-graded $K$-algebra $A$ we have
identified the category $\bLMod_A$ with the category of left actions of $A$.
A left $A$ module $N$ corresponding to the abelian object in $\bCA_K/A$ 
given by
$\pr_1:A\oplus N\downarrow A$ with unit section $a\mapsto(a,0)$ and
product given by $(a,x)(b,y)=(ab,ay+xb)$. A section of this object of
$\bCA_K/A$ is a (degree-preserving) derivation, that is, an $X$-graded
$K$-module map $d:A\to N$ such that 
$d(ab)=a\,db+b\,da$, as usual. The functor $N\mapsto\Der_K(A,N)$ is 
co-represented by the (graded) $A$-module of 
K\"ahler differentials: $\Omega_{A/K}\in\bLMod_A$. 
Expressed in terms of generators and relations, this $A$-module is the
cokernel of the map 
\[
d:\bigoplus_{x,y}P^{x+y}\to\bigoplus_zP^z
\]
determined by 
\[
d\circ\inj_{x,y}=\inj_x\circ y^*-\inj_{x+y}+\inj_y\circ x^*\,.
\]
This is $\Ab_AA$. To describe $\Ab_AB$, for $p:B\to A$ in $\bCA(\bMod_K^X)$, 
notice that for each $x\in X$ the right $A$-module $P_x$ can be 
regarded as a right $B$-module through the map $p$. Then 
\[
(\Ab_AB)_x=P_x\tensor_B\Omega_{B/K}\,.
\]
The left $A$-module structure on $\Ab_AB$ arises from the $A$-bimodule
structure of $P$.

The $(F_X,u)$ adjoint pair of Section \ref{sec-gradings} produces
a cotriple on $\bCA(\bMod_K^X)$ which we denote by $\Sym^X_K$.  
Following \cite{barr-beck}, 
this in turn leads to a natural simplicial object $\Sym^X_{K\bullet}A$
augmented to $A$, with $\Sym^X_{Kn}A=(\Sym^X_K)^{n+1}A$,
the simplicial cotriple resolution, which can be used to derive 
functors on $\bCA(\bMod_K^X)/A$. 
The Quillen homology of $A\in\bCA(\bMod^X_K)$ is defined as the derived
functors of $\Ab_A:\bCA(\bMod^X_K)/A\to\bLMod_A$, evaluated at the object
$\id_A:A\downarrow A$.

For each $n\geq0$ the Quillen homology $HQ_n(A)$ is itself an $A$-module,
and $HQ_0(A)=\Omega_{A/K}$. We can endow the Quillen homology with 
coefficients in a right $A$-module $M$: 
\[
HQ_n(A;M)=\pi_n(M\tensor_A\Ab_A\Sym^X_{K\bullet}A)\,.
\]
For any $x\in X$, we can recover
$HQ_n(A)_x$ by using the right $A$-module $P_x$ for coefficients:
\[
HQ_n(A)_x=HQ_n(A;P_x)\,.
\]
For $N\in\bLMod_A$  we can define the Andr\'e-Quillen cohomology as follows:
\[
HQ^n(A;N)=H^n(\Hom_A(\ch\Sym^X_{K\bullet}A,N))\,.
\]

\section{Change of grading monoid}

Let $\alpha:Y\to X$ be a map of commutative monoids. An $X$-graded object 
$C$ in $\bC$  determines a $Y$-graded object $\alpha^*C$ by
\[
(\alpha^*C)_y=C_{\alpha(y)}\,.
\]
If $\bC$ has coproducts, 
the functor $\alpha^*:\bC^X\to\bC^Y$ has a left adjoint given by
\[
(\alpha_*C)_x=\coprod_{\alpha(y)=x}C_y\,.
\]

Let $K$ be a commutative ring.
The functor $\alpha_*:\bMod_K^Y\to\bMod_K^X$ is symmetric monoidal --
\[
\alpha_*(M\tensor N)=\alpha_*M\tensor\alpha_*N
\]
-- so $\alpha_*$ sends $Y$-graded commutative $K$-algebras to $X$-graded
commutative $K$-algebras. Then for any $A\in\bCA(\bMod^Y_K)$ the map $\alpha$ 
induces adjoint pairs
\[
\alpha_*:\bLMod_A\rightleftarrows\bLMod_{\alpha_*A}:\alpha^*
\]
and
\[
\alpha_*:\bRMod_A\rightleftarrows\bRMod_{\alpha_*A}:\alpha^*
\]

Let $N$ be an $\alpha_*A$-module. It is straightforward to construct a natural
isomorphism
\[
\Der_K(A,\alpha^*N)=\Der_K(\alpha_*A,N)
\]
and thus a natural isomorphism of $\alpha_*A$-modules
\[
\Omega_{\alpha_*A/K}=\alpha_*\Omega_{A/K}\,.
\]

The adjoint pairs $(F_X,u_X)$ are compatible under change of grading monoid:
A monoid homomorphism $\alpha:Y\to X$ determines the following squares.
\[
\xymatrix{
\bSet^Y & \bCA(\bMod^Y_K) \ar[l]_-{u_X} &
\bSet^Y \ar[r]^-{F_Y} \ar[d]^{\alpha_*} & \bCA(\bMod^Y_K) \ar[d]^{\alpha_*}  \\
\bSet^X \ar[u]_{\alpha^*} & \bCA(\bMod^X_K) \ar[l]_-{u_X} \ar[u]_{\alpha^*} &
\bSet^X \ar[r]^-{F_X}  & \bCA(\bMod^X_K) 
}
\]
The diagram of right adjoints clearly commutes, so the diagram of left 
adjoints does too:
\[
\alpha_*F_Y(T)=F_X(\alpha_*T)\in\bCA(\bMod^X_K)\,.
\]
Passing to the cotriples,
\[
\alpha_*\Sym_K^Y(A)=\alpha_*F_Yu_X(A)=F_X\alpha_*u_Y(A)=F_Xu_X(\alpha_*A)=
\Sym_K^X(\alpha_*A)
\]
and thus to simplicial resolutions:
\[
\alpha_*\Sym^Y_{K\bullet}(A)=\Sym^X_{K\bullet}(\alpha_*A)\,.
\]

Assembling all this, along with the fact that $\alpha_*$ is exact, we find:

\begin{proposition} Let $\alpha:Y\to X$ be a homomorphism of commutative 
monoids. 
There are isomorphisms of $\alpha_*A$-modules natural in $A\in\bCA(\bMod_K^Y)$
\[
\alpha_*HQ_n(A)=HQ_n(\alpha_*A)
\]
as well as isomorphisms of $K$-modules 
\[
HQ_n(A;\alpha^*M)=HQ_n(\alpha_*A;M)
\]
natural in $M\in\bRMod_{\alpha_*A}$ and 
\[
HQ^n(A;\alpha^*N)=HQ^n(\alpha_*A;N)
\]
natural in $N\in\bLMod_{\alpha_*A}$.
\end{proposition}

An important example is provided by taking $X$ to be the one-element monoid,
$e$, and $\alpha:Y\to e$ the unique map.  Then $\alpha_*A$ is the ``degrading''
of $A$, an ungraded commutative $K$-algebra, and $N$ is a module for it;
$\alpha^*N$ is the ``constant'' $Y$-graded $K$-module with 
$(\alpha^*N)_y=N$ for all $y\in Y$, and $A$ acting among them in the obvious
way; and $HQ_*(\alpha_*A)$, $HQ_*(\alpha_*A;M)$, and $HQ^*(\alpha_*A;N)$ 
are the usual Andr\'e-Quillen groups. 

\section{Commutative monoids} \label{sec-cm}
We now regard commutative monoids as the category of homological interest,
rather than as a source of gradings. 

Let $X\in\bCM$. It is easy \cite{leech,grillet-74} 
to identify the category of Beck modules over $X$ 
in terms of the {\em Leech category}, $L_X$, 
with object set $X$ and $L_X(x,z)=\{y\in X:x+y=z\}$ with unit and 
composition determined by the commutative monoid structure on $X$. Write 
$y_*:x\to(x+y)$ for a morphism in this category. The category of Beck 
modules over $X$ is canonically equivalent to the category of functors
from $L_X$ to the category $\bAb$ of abelian groups. 

A map $\alpha:Y\to X$ of commutative monoids induces a functor
\[
\alpha^*:\bLMod_X\to\bLMod_Y
\]
which, under the equivalence with functors from the Leech categories, 
may be regarded as induced by pre-composition with the induced functor 
$\alpha:L_Y\to L_X$. The left adjoint $\alpha_*$ is then induced by 
left Kan extension \cite[Chapter X]{maclane} along $\alpha$. 

A section of an abelian object over $X$ is a ``derivation,'' and under
the identification of abelian objects over $X$ with functors on $L_X$ a 
derivation with values in $N:L_X\to\bAb$ is an assignment of an element
$s(x)\in N_x$ for each $x\in X$ such that
\[
s(x+y)=x_*s(y)+y_*s(x)
\]
There is a universal example, the Beck module of ``K\"ahler differentials''
$\Omega^{CM}_X$, which provides a distinguished object of $\bLMod_X$. 
For any $\alpha:Y\to X$, 
\[
\Ab_XY=\alpha_*\Omega^{CM}_Y\,.
\]

A commutative monoid $X$ defines a canonical commutative $X$-graded $K$-algebra
$\wK X$ in which $(\wK X)_x=K$ for each $x\in X$, with generator $1_x$,
$1=1_0\in(\wK X)_0$, and 
\[
\mu_{x,y}(a_x\tensor b_y)=(ab)_{x+y}=ab1_{x+y}\,,\quad a,b\in K\,.
\]
This object of $\bCA_K^X$ co-represents the functor sending an object $A$
to the set of sections $a$ of the grading 
function $A\to X$ such that $a_xa_y=a_{x+y}$.

For a Beck $K$-module $N$ over $X$ we may define a left module over 
$\wK X$, which we will also denote by $N$. The $K$-module in degree $x$
is $N_x$ and the multiplication by $\wK X$ is given by the action of elements
of $X$, regarded as morphisms in $L_X$. The following lemma is simply a change
in perspective.

\begin{lemma} 
The category of covariant functors $L_X\to\bMod_K$ is canonically
equivalent to the category of $\wK X$-modules. The category of 
contravariant functors $L_X\to\bMod_K$ is canonically equivalent to the
category of right $\wK X$-modules.
\end{lemma}

Of course, the usual monoid-algebra is obtained from $\wK X$ as
\[
KX=\alpha_*\wK X\,,\quad\alpha:X\to e\,.
\]
These isomorphisms are compatible with attendant structure: the functors
induced by a map $\alpha:Y\to X$ of commutative monoids agree with the 
functors induced by $\wZZ\alpha:\wZZ Y\to\wZZ X$. Under this identification
the K\"ahler differential objects match up:
\[
\Omega^{CM}_Y=\Omega^{CA}_{\wZZ Y}
\]
and more generally $\alpha_*\wZZ Y\in\bCA(\bMod^X_K)/\wZZ X$ and 
\[
\Ab_X(Y)=\Ab_{\wZZ X}(\alpha_*\wZZ Y)
\]
in $\bLMod_X=\bLMod_{\wZZ X}$.
As a further example, given a contravariant functor $M:L_X\to\bMod_K$ and a 
covariant functor $N:L_X\to\bMod_K$,
\[
M\tensor_{\wZZ X}N=M\tensor_{L_X}N
\]
where the right hand side is the usual tensor product over the category $L_X$,
as considered by Kurdiani and Pirashvili in \cite{kurdiani-pirashvili}.

In order to express the functoriality of $\wK$,
we need a category of graded algebras in which the grading commutative
monoid can vary. So let $\bCA_K^*$ be the category whose objects are
pairs $(X,A)$ where $X$ is a commutative monoid and $A$ is a commutative 
$X$-graded $K$-algebra. A morphism $(X,A)\to(Y,B)$ consists of a monoid
homomorphism $\alpha:X\to Y$ together with a morphism
$f:\alpha_*A\to B$ in $\bCA(\bMod_K^Y)$ (or equivalently a morphism
$\hat f:A\to\alpha^*B$ in $\bCA(\bMod_K^X)$). Given a second morphism
$(\beta,g):(Y,B)\to(Z,C)$, the composite is defined as 
$(\beta\circ\alpha,g\circ\beta_*f)$, (or equivalently as
$(\beta\circ\alpha,\alpha^*\hat g\circ\hat f)$). Then the construction
$\wK$ provides a functor
\[
\wK:\bCM\to\bCA_K^*\,.
\]

Here is the main theorem of the paper. 
\begin{theorem} \label{cm-ca}
There are isomorphisms, natural in the commutative monoid $X$,
\[
HQ^{CM}_n(X)=HQ^{CA}_n(\wZZ X)\,,
\]
where we have identified Beck modules over $X$ with modules over $\wZZ X$. 
For any right Beck module $M$ over $X$,
\[
HQ^{CM}_n(X;M)=HQ^{CA}_n(\wZZ X;M)
\]
and for any left Beck module $N$
\[
HQ_{CM}^n(X;N)=HQ_{CA}^n(\wZZ X;N)\,.
\]
\end{theorem}
\begin{proof}
Let $\alpha:Y_\bullet\to X$ be a cofibrant replacement of $X\in\bCM$ 
regarded as a constant simplicial object. 
We could take the cotriple resolution $\NN_\bullet X$ for example,
and for definiteness we will do so. We claim that the induced map
$\alpha_*\wK\NN_\bullet X\to\wK X$ is a cofibrant replacement for $\wK X$
in $s\bCA(\bMod_K^X)$. 

The augmentation $\alpha:\NN_\bullet X\to X$
is a weak equivalence in $s\bCM$, which is to say a weak
equivalence of simplicial sets. The entire simplicial set
$\NN_\bullet X$ splits into a disjoint union of the pre-images under $\alpha$
of the elements of $X$, 
and each of these components is a weakly contractible simplicial set. 

Now apply the functor $\wK$ to $\NN_\bullet X$, to get a simplicial object in 
$\bCA_K^*$ with an augmentation $\wK\NN_\bullet X\to\wK X$; that is, a map 
\[
\alpha_*\wK\NN_\bullet X\to\wK X
\]
in $s\bCA(\bMod^X_K)$. The object $\alpha_*\wK\NN_\bullet X$ of $s\bMod^X_K$ 
splits as a direct sum of the free $K$-module functor applied 
to the components of $\NN_\bullet X$ above the elements of $X$. 
But a weak equivalence $Z_\bullet\to Y_\bullet$ of simplicial sets
induces a weak equivalence $KZ_\bullet\to KY_\bullet$
\cite[III, Prop. 2.16]{goerss-jardine}, so
the map $\alpha_*\wK\NN_\bullet X\to\wK X$ is a weak equivalence.

The second observation is that $\alpha_*\wK\NN_\bullet X$ is cofibrant
in the model category of simplicial objects in $\bCA(\bMod_K^X)$. 
Indeed, it is almost free in the sense of \cite{miller}, and hence,
as observed there, is cofibrant. 
A simplicial object $Z_\bullet\in s\bCM$ is ``almost free'' if there are
subsets $G_s\subseteq Z_s$, for each $s$, that are respected by all 
face and degeneracy maps except for $d_0$, and such that $Z_s$ is 
freely generated by $G_s$. A cotriple resolution, such as $\NN_\bullet X$,
is easily seen to be almost free. If we then apply the functor 
$\alpha_*\wK$ to it, the same generating sets
again generate as objects in $\bCA(\bMod^X_K)$; it is 
an almost free object in $s\bCA(\bMod_K^X)$.

The map $\alpha_*\wK\NN_\bullet X\to\wK X$ thus joins
$\Sym^X_{K\bullet}(\wK X)\to\wK X$ as a cofibrant replacement, and
so (e.g. \cite[Proposition 3.9]{frankland}) can be used to compute
Quillen homology. Thus (finally taking $K=\ZZ$)
\[
\Ab_X\NN_\bullet X=\Ab_{\wZZ X}(\alpha_*\wZZ\NN_\bullet X)
\simeq\Ab_{\wZZ X}(\Sym_{\ZZ\bullet}^X(\wZZ X))
\]
as simplicial objects in $\bLMod_X=\bLMod_{\wZZ X}$.
Applying the functors $\pi_*(-)=H_*(\ch(-))$,
$\pi_*(M\tensor_{\wZZ X}-)=H_*(M\tensor_{\wZZ X}\ch(-))$,
and $H^*(\Hom_{\wZZ X}(\ch(-),N))$ gives us the results. 
\end{proof}

\begin{corollary}[\cite{kurdiani-pirashvili}]
Let $X$ be a commutative monoid and 
let $\alpha:X\to e$ be the unique monoid map to the trivial monoid. 
There are isomorphisms
\[
\alpha_*H^{CM}_n(X;\alpha^* M)=H^{CA}_n(\ZZ X;M)
\]
natural in the right $X$-module $M$, and
\[
\alpha_*H_{CM}^n(X;\alpha^* N)=H_{CA}^n(\ZZ X;N)
\]
natural in the left $X$-module $N$.
\end{corollary}

\section{The Hochschild complex}
This section reviews well known material (e.g. \cite{weibel,witherspoon}) 
in the graded setting. 

Let $A$ be an associative $X$-graded $K$-algebra. There is a canonical
simplicial $A$-bimodule (so each level is an $X$-graded $K$-module with
grade-preserving left and right actions of $A$) $B_\bullet(A)$ over $K$ with 
\[
B_n(A)=A^{\tensor(n+2)}
\]
for $n\geq0$, augmented to $A$, and 
\[
d_i=1^{\tensor i}\tensor \mu\tensor1^{\tensor(n-i)}:A^{\tensor(n+2)}
\rightarrow A^{\tensor(n+1)}\,,\,\, 0\leq i\leq n
\]
\[ 
s_i=1^{\tensor(i+1)}\tensor\eta\tensor1^{\tensor(n-i)}:A^{\tensor(n+1)}
\rightarrow A^{\tensor(n+2)}\,,\,\, 0\leq i\leq n-1
\]
where $\mu:A\tensor A\rightarrow A$ is the multiplication and 
$\eta:K\rightarrow A$ includes the unit. 

There are also maps 
\[
s_{-1}=\eta\tensor1^{\tensor(n+1)}\,,\quad s_n=1^{\tensor(n+1)}\tensor\eta:
A^{\tensor(n+1)}\rightarrow A^{\tensor(n+2)}\,.
\]
The first is a right $A$-module map and the second is a left $A$-module map,
and they provide contracting homotopies of the simplicial object regarded
as either a simplicial right $A$-module or a simplicial left $A$-module.
In fact it is just the simplicial
bar resolution of $A$ as a left or right $A$-module. 
Thus the chain complex associated to $B_\bullet(A)$, $\ch B_\bullet(A)$, 
is a relative projective resolution of $A$ as an 
$A$-bimodule, the {\em Hochschild resolution}. 
If $A$ is projective as a $K$-module, it is an absolute
projective resolution. 

Let $Q_A$ be the functor from $A$-bimodules to $X$-graded $K$-modules with
\[
Q_A(M)_z=M/K\{am-ma:a\in A_x,m\in M_w,w+x=z\}\,.
\]
More generally, given an $A$-bimodule $N$, define the functor $Q^A_N$, or
$Q_N$ if $A$ is understood, from $A$-bimodules to $X$-graded $K$-modules
\[
Q_N(M)_z=(M\tensor N)_z/R
\]
where $R$ is the sub-$K$-module generated by
\[
\{am\tensor n-m\tensor na,ma\tensor n-m\tensor an:
m\in M_w, a\in A_x, n\in N_y,w+x+y=z\}\,.
\]
We recover $Q_A$ by regarding $A$ as a bimodule over itself using left
and right multiplication. 
A bimodule is the same thing as a module over 
$A^e=A\tensor A^{op}$, and under this equivalence
\[
Q_N(M)=M\tensor_{A^e}N\,.
\]
In general, this is just an $X$-graded $K$-module, but if $A$ is commutative 
then $Q_A(M)$
is naturally an $A$-module, since $A$ is then an $(A,A^e)$-bimodule.

Apply this functor to $B_\bullet(A)$ to obtain a simplicial object 
in $\bMod^X_K$, $Q_NB_\bullet(A)$, equipped with an augmentation to $Q_NA$.
This is the {\em Hochschild complex} with coefficients in $N$,
\[
C_\bullet(A/K;N)=\ch Q_NB_\bullet(A)
\]
and by definition
\[
Hoch_n(A/K;N)=H_n(C_\bullet(A/K;N))\,.
\]
When the ground ring $K$ is understood we may drop it from the notation.
When $N=A$ with its natural $A$-bimodule structure, we may drop it from the
notation as well:
\[
C_\bullet(A)=C_\bullet(A/K;A)\quad\text{and}\quad
Hoch_\bullet(A)=Hoch_\bullet(A/K;A)\,.
\]

To understand this better, notice the isomorphism
\[
Q_N(A\tensor V\tensor A)\rightarrow V\tensor N
\]
for $V\in\bMod^X_K$, given by factoring 
\[
a\tensor v\tensor b\tensor n\mapsto v\tensor bna
\]
through $Q_N(A\tensor V\tensor A)$. The inverse sends
$v\tensor n$ to $[1\tensor v\tensor1\tensor n]$.

This isomorphism breaks symmetry. But using it we may write the augmented
simplicial $K$-module  $Q_NB_\bullet(A)$ as
\[
Q_N(A)\leftarrow N\Leftarrow A\tensor_KN\Lleftarrow\cdots\,,
\]
so
\[
C_n(A;N)=A^{\tensor n}\tensor_KN\,.
\]

If $A$ and $Z$ are two $X$-graded $K$-algebras and $M$ and $N$ bimodules 
for them, there is a natural isomorphism
\[
Q_A(M)\tensor Q_Z(N)\rightarrow Q_{A\tensor Z}(M\tensor N)
\]
under the identity map on $M\tensor N$. The fact that it is an isomorphism
follows from the identity
\[
(a\tensor z)(m\tensor n)-(m\tensor n)(a\tensor z)
=am\tensor(zn-nz)+(am-ma)\tensor nz
\]
We get natural isomorphisms of simplicial objects
\begin{gather*}
B_\bullet(A)\tensor B_\bullet(Z)\rightarrow B_\bullet(A\tensor Z)\\
Q_AB_\bullet(A)\tensor Q_ZB_\bullet(Z)\rightarrow 
Q_{A\tensor Z}B_\bullet(A\tensor Z)
\end{gather*}

If $A$ is commutative, we may take $A=Z$ and compose with the $K$-algebra map
$\mu:A\tensor A\rightarrow A$ to obtain a simplicial commutative $A$-algebra 
structure on $Q_AB_\bullet(A)$, and $Q_NB_\bullet(A)$
becomes a module over $Q_AB_\bullet(A)$. 

Passing to associated chain complexes, the Eilenberg-Zilber or shuffle map 
(\cite[p.~64]{eilenberg-maclane} or \cite[p.~39]{miller-lectures})
results in the structure of a commutative (in the signed sense) 
differential graded $A$-algebra 
on $C_\bullet(A)$ and hence a graded commutative $A$-algebra 
structure on its homology $Hoch_\bullet(A)$.  

Dually, the Hochschild cochain complex with coefficients in 
and $A$-bimodule $N$ is  $\Hom_{A\tensor A^{op}}(B_\bullet(A),N)$, 
and its homology is the Hochschild cohomology $Hoch^\bullet(A;N)$. 

It is well known and easy to verify that
\[
Hoch_1(A)=\Omega_{A/K}=HQ^{CA}_0(A)
\]
and 
\[
Hoch^1(A;M)=\Der_K(A;M)=HQ_{CA}^0(A;M)\,.
\]

\section{Harrison homology and cohomology}

Now suppose that $A$ is a commutative $K$-algebra. 
Then $Q_AA=A$; the Hochschild complex $C_\bullet(A)$
is augmented to $A$. Let $I_\bullet(A)$ denote the kernel of this augmentation;
this is the ideal of positive-dimensional elements 
in the commutative differential graded $A$-algebra $C_\bullet(A)$.
The {\em Harrison complex} \cite{harrison} 
is the differential graded module of 
indecomposables in $C_\bullet(A)$, $I_\bullet(A)/I_\bullet(A)^2$.
The {\em Harrison homology} of $A$ is the homology of this chain complex
of  $A$-modules:
\[
Harr_n(A)=H_n(I_\bullet(A)/I_\bullet(A)^2)\,.
\]
We can equip it with coefficients in an $A$-module $M$:
\[
Harr_n(A;M)=H_n((I_\bullet(A)/I_\bullet(A)^2)\tensor_AM)\,.
\]
The Harrison cohomology with coefficients in an $A$-module $M$ is 
\[
Harr^n(A;M)=H^n(\Hom_A(I_\bullet(A)/I_\bullet(A)^2,M))\,.
\]
Clearly 
\[
Harr_0(A)=0\quad\text{and}\quad Harr_1(A)=Hoch_1(A)
\]

The shuffle product defines a sign-commutative graded $K$-algebra structure on 
the $\NN$-graded $K$-module $\oC_\bullet(A)$ with
\[
\oC_n(A)=A^{\tensor n}\,.
\]
As graded $A$-algebras
\[
C_\bullet(A)=A\tensor\oC_\bullet(A)\,.
\]
Only the differential depends on the algebra
structure, and it is not the $A$-linear
extension of a differential on $\oC_\bullet(A)$. 

The $K$-module of Harrison cochains can be re-expressed in terms of the 
graded $K$-algebra
$\oC_\bullet(A)$. 
Let $\oI_\bullet(A)$ be its augmentation ideal; then
\[
I_\bullet=A\tensor\oI_\bullet\,,\quad I_\bullet^2=A\tensor\oI_\bullet^2\,,
\]
and so
\[
I_\bullet(A)/I_\bullet(A)^2=A\tensor(\oI_\bullet(A)/\oI_\bullet(A)^2)\,.
\]
A Harrison $n$-cochain (for $n>0$) with coefficients in $M$ 
is thus a $K$-linear map
\[
s:A^{\tensor n}\to M
\]
that annihilates decomposables. This may be phrased
as a symmetry condition on the cochain: given $i,j$, both positive and
summing to $n$, let $\Sigma(i,j)$ be the set $(i,j)$-shuffles; that is, the
set of elements of $\Sigma_n$ 
that preserve the order of $\{1,\ldots,i\}$ and of $\{i+1,\ldots,n\}$. 
The symmetry condition $(i,j)$ on a Hochschild cochain $s:A^{\tensor n}\to M$
is 
\[
\sum_{\sigma\in\Sigma(i,j)}\sgn(\sigma)s\circ\sigma = 0\,.
\]
Since the shuffle product is commutative, we may assume that $i\leq j$;
there are $\lfloor n/2\rfloor$ independent conditions.  

An alternative symmetry condition (apparently the one originally conceived 
of by Harrison; the shuffle description is said to be due to Mac Lane) 
is described in \cite{gerstenhaber-schack}. Think of an element of $\Sigma_n$
as an ordering of $\{1,2,\ldots,n\}$. Let $1\leq k\leq n$. An 
element $\sigma\in\Sigma_n$ is a $k${\em -monotone} permutation if the lead
element is $k$, the numbers $1,2,\ldots,k$ occur in decreasing order,
and the numbers $k+1,\ldots,n$ occur in increasing order.
There are $\binom{n-1}{k-1}$ of them. 
For example there is only one $n$-monotone permutation in $\Sigma_n$, 
corresponding to the
sequence $n,n-1,\cdots,2,1$, and the $4$-monotone permutations in $\Sigma_6$ 
are
\[
432156, 432516, 432561, 435216, 435261, 435621, 453216, 453261, 453621, 456321
\]
Let $M_k(n)$ be the set of $k$-monotone permutations in $\Sigma_n$.
Let $dr(\sigma)$ be the sum of the positions occupied by $1,2,\ldots,k-1$
in the permutation $\sigma\in M_k(n)$.

\begin{lemma}[\cite{gerstenhaber-schack}, Theorem 4.1] 
A map $s:A^{\tensor n}\to M$ is a Harrison cochain if and only if 
\[
s=\sum_{\sigma\in M_k(n)}(-1)^{dr(\sigma)}s\circ\sigma\,,\quad2\leq k\leq n\,.
\]
\end{lemma}
So for example a Hochschild $4$-cochain $s$ is a Harrison cochain 
if and only if
\begin{align*}
s(a_1,a_2,a_3,a_4)&=s(a_2,a_1,a_3,a_4)-s(a_2,a_3,a_1,a_4)+s(a_2,a_3,a_4,a_1)\\
&=-s(a_3,a_2,a_1,a_4)+s(a_3,a_2,a_4,a_1)-s(a_3,a_4,a_2,a_1)\\
&=-s(a_4,a_3,a_2,a_1)
\end{align*}
The $4$-monotone and $2$-monotone symmetries combine to give
\[
s(a_4,a_3,a_2,a_1)=-s(a_2,a_1,a_3,a_4)+s(a_2,a_3,a_1,a_4)-s(a_2,a_3,a_4,a_1)
\]
which is the same as the $3$-monotone condition (after rearranging the labels);
so we can dispense with either one of the first two conditions in this list.
The same argument shows that one need only assume the $n$-monotone
condition together with one condition from each pair
$\{2,n-1\},\{3,n-2\},\ldots$: so $\lfloor n/2\rfloor$ conditions suffice. 
This matches with the number of independent shuffle conditions.

\section{Divided powers and Barr homology}

Michael Barr suggested possible variations on Harrison's symmetry conditions, 
in an attempt to come closer to Quillen homology. As explained by 
Sarah Whitehouse \cite{whitehouse}, these variations still fail, though
they may give better approximations. 

Gerstenhaber and Schack 
\cite[Remark, p.~232]{gerstenhaber-schack} (see also \cite{whitehouse})
suggest that one of Barr's ideas was to 
divide the Hochschild complex not just by shuffle decomposables but by 
the divided power structure as well. While a divided power structure
on the even homotopy groups of a simplicial commutative algebra was implicit 
in the works of Eilenberg and Mac Lane \cite{eilenberg-maclane} and 
Henri Cartan \cite[Exp. 8]{cartan}, 
its construction on the level of the
Hochschild complex was at best a folk result at the time of Barr's question,
and even when Gerstenhaber and Schack were writing. It seems to have
first been set out, in the associated chain complex of a simplicial
commutative algebra $B_\bullet$, by Siegfried Br\"uderle and Ernst Kunz 
\cite{bruederle-kunz} in 1994; see also \cite{richter} and \cite{gillam}. 
The result is a natural family of maps
\[
\gamma_k:B_{2n}\to B_{2kn}
\]
such that $k!\gamma_kx=x^k$.

From these sources one obtains the following formula for the divided power
structure on $\overline C_{even}(A)$, where $A$ is a commutative $K$-algebra.
Let $S_k(kn)$ be the set of shuffles associated to the partition of 
$\{1,2,\ldots,kn\}$ into $k$ intervals of length $n$.
Let $S'_k(kn)$ be the subset of these
such that the leading terms of the $k$ sequences occur in order. 
Then, for $n$ even,
\[
\gamma_k[a_1|\cdots|a_n]=\sum_{\sigma\in S'_k(kn)}\sgn(\sigma)
[a_1|\cdots|a_n|a_1|\cdots|a_n|\cdots|a_1|\cdots|a_n]\circ\sigma
\]
where the sequence $a_1|\cdots|a_n$ is repeated $k$ times. For example
\begin{gather*}
\gamma_2[a|b]=[a|b|a|b]\,,\quad \gamma_3[a|b]=[a|b|a|b|a|b]\,, \\
\gamma_2[a|b|c|d]=[a|b|c|d|a|b|c|d]-[a|b|c|a|d|b|c|d]+[a|b|c|a|b|d|c|d] \\
+[a|b|a|c|d|b|c|d]-[a|b|a|c|b|d|c|d]+2[a|b|a|b|c|d|c|d]\,.
\end{gather*}

We can put at least one restriction on the tensors occurring in the expression 
for the divided powers. To express it, notice that there is a universal
Hochschild $n$-chain, $[a_1|\cdots|a_n]\in K[a_1,\ldots,a_n]^{\tensor n}$.

\begin{lemma} No decomposable tensor with entries chosen from 
$\{a_1,\ldots,a_n\}$ occurring with nonzero coefficient in
$\gamma_k[a_1|\cdots|a_n]$ has consecutive occurrences of any $a_i$.
\end{lemma}
\begin{proof} We show how such terms cancel in pairs in the expression for the
divided power, by defining a free involution on the set of terms with
neighboring repeated letters such that the elements of each
orbit occur with opposite signs. The involution will leave unchanged
all the letters up to and including the left-most neighboring repeated pair.

If the repeated pair is $a_1|a_1$, swap the positions of the remaining
letters in the two blocks initiated by these letters. We get an identical
word, but since $n$ is even
this is an odd number of transpositions, so the terms cancel. 

If the repeated pair is $a_i|a_i$ for $i>1$, just swap those two entries.
This is allowed since the leading term of both blocks precedes both entries
in the repeated pair.
\end{proof}

Since every term in the expression for $\gamma_k[a_1|\cdots|a_n]$ has 
leading entry $a_1$, we obtain:

\begin{corollary} $\gamma_k[a|b]=[a|b|a|b|\cdots|a|b]$.
\end{corollary}

In general the expression for $\gamma_k$ seems very complicated. 
For example, a computer calculation shows that 
$\gamma_3[a|b|c|d]$ has 53 terms, with coefficients ranging from $-4$ to $6$.

Write $Barr_*(A)$ for the homology of the chain complex of $A$-modules
obtained from $C_\bullet(A)$ quotienting out by decomposables and the image
of divided powers.
Since $k!\gamma_k(\omega)$ is decomposable, we have an exact sequence
\[
0\to Harr_{2k+1}(A)\to Barr_{2k+1}(A)\to T_{2k}\to Harr_{2k}(A)
\to Barr_{2k}(A)\to0
\]
where $k!T_{2k}=0$. Thus
\begin{gather*}
Harr_{2k+1}(A)\to Barr_{2k+1}(A)\,\,
\text{is injective with cokernel killed by $k!$}\,,\\
Harr_{2k}(A)\to Barr_{2k}(A)\,\,
\text{is surjective with kernel killed by $k!$}\,.
\end{gather*}

We can also form the ``Barr cohomology'' with coefficients in an $A$-module. 
Its cochains are the 
Hochschild cochains $s$ satisfying the Harrison symmetry conditions
with the additional conditions
\[
s(\gamma_k(\omega))=0\,,\quad k>0\,,
\]
for $|\omega|$ even; for example $s(a,b,a,b)=0$ in dimension 4;  
$s(a,b,a,b,a,b)=0$ in dimension 6; and in dimension 8 
there are two additional symmetries, 
\[
s(a,b,a,b,a,b,a,b)=0
\]
guaranteeing annihilation of $\gamma_4[a|b]$, and
\begin{gather*}
s(a,b,c,d,a,b,c,d)-s(a,b,c,a,d,b,c,d)+s(a,b,c,a,b,d,c,d)\\
+s(a,b,a,c,d,b,c,d)-s(a,b,a,c,b,d,c,d)+2s(a,b,a,b,c,d,c,d)=0
\end{gather*}
to annihilate $\gamma_2[a|b|c|d]$. 

\begin{remark}{\em
It is natural to hope that the natural map of 
$A$-modules $Hoch_n(A)\to HQ_{n-1}(A)$ factors as
\[
Hoch_n(A)\to Harr_n(A)\to Barr_n(A)\to HQ_{n-1}(A)\,,
\]
but this seems unlikely to us except in low dimensions.
}\end{remark}

\section{Grillet's work} \label{sec-grillet}

In a series of papers, Pierre Grillet associated to a commutative monoid $X$
and a Beck module $M$ over it the beginning of a cochain 
complex and proves or conjectures that it computes the low-dimensional 
components of the Quillen cohomology
$HQ_{CM}^*(X;M)$. We observe that the symmetry conditions he imposed are
precisely the monotone conditions, with two variations corresponding to
Barr's variation on the Harrison complex. 
We will not attempt a complete survey of Grillet's work on this subject,
but merely note the occurrence of symmetry conditions that we now
see as Harrison or Barr symmetry conditions on 
Hochschild cochains, and where in Grillet's work
they are proved to yield cohomology groups isomorphic to Quillen's.

To begin with, for any $X$ graded $K$-algebra and $A$-module $M$,
\[
Harr^1(A;M)=Hoch^1(A;M)=HQ_{CA}^0(A;M)= Der_K(A,M)
\]
so
\[
Harr^1(\wZZ X;M)= HQ_{CM}^0(X;M)\,.
\]

In 1974 \cite{grillet-74} Grillet used 2-cocycles $s$ with the symmetry
\[
s(a,b)=s(b,a)
\]
to classify extensions of commutative monoids. This is of course the 
$2$-monotone symmetry. Twenty years later, in 
\cite{grillet-95}, he returned to this by invoking Quillen cohomology
as an intermediary, thus showing that 
\[
Harr^2(\wZZ X;M)=HQ_{CM}^1(X;M)\,.
\]
(Grillet chooses to index Quillen homology following the Hochschild convention,
so he would write $HQ^2_{CM}(X;M)$.)

This study was supplemented by \cite{grillet-97}, which confirmed this
result by direct computation and extended it to dimension 3 using the
symmetry conditions
\begin{gather*}
s(a,b,c)+s(c,b,a)=0\\
s(a,b,c)+s(b,c,a)+s(c,a,b)=0\,.
\end{gather*}
Taken together, these are equivalent to the $k$-monotone symmetries for 
$k=2$ and $3$. With a remarkable calculation Grillet then verified that
\[
Harr^3(\wZZ X;M)=HQ^2_{CM}(X;M)
\]
This work was consolidated and summarized in his book \cite{grillet-01}.

After another twenty years, Grillet returned again to this project, 
in \cite{grillet-20}, extending his calculation to dimension 4 using 
cochains satisfying the symmetry conditions simplified in 
\cite{grillet-21} to
\begin{gather*}
s(a,b,c,d)-s(b,a,c,d)+s(b,c,a,d)-s(b,d,d,a)=0\,,\\
s(a,b,c,d)+s(d,c,b,a)=0\,,\\
s(a,b,b,a)=0\,.
\end{gather*}
The reader will recognize the first two as the $2$-monotone and $4$-monotone
symmetries.
The first equation implies $s(a,b,a,b)=s(b,a,a,b)$, so the third condition
is equivalent to the Barr variant $s(a,b,a,b)=0$.
Here again, Grillet obtained the surprising result
\[
Barr^4(\wZZ X;M)=HQ^3_{CM}(X;M)\,.
\]

This work was quickly followed by \cite{grillet-22}, in which Grillet 
proposed symmetric conditions on Hochschild cocycles extending into
dimensions 5 and 6. In dimension 5, he proposed
\begin{gather*}
s(a,b,c,d,e)-s(b,a,c,d,e)+s(b,c,a,d,e)-s(b,c,d,a,e)+s(b,c,d,e,a)=0\\
s(a,b,c,d,e)+s(e,d,c,b,a)=0\,.
\end{gather*} 
We recognize these as the $2$-monotone and $5$-monotone conditions, 
which suffice to determine Harrison cohomology. In dimension 6, his 
proposed symmetries are precisely the $k$-monotone conditions for $k=2,3$, 
and 6, augmented by the Barr variant $s(a,b,a,b,a,b)=0$. 

In these higher dimensions
the identifications with Quillen cohomology are left as conjectures.
We can now see that at least the 5-dimensional case was too optimistic.
Write $\alpha$ for the unique map of commutative monoids $\NN\to e$. Then
\[
Barr^\bullet(\wZZ\NN;\alpha^*\FF_p)=Barr^\bullet(\ZZ[x];\FF_p)\,.
\]
Let $c:\FF_2[x]^{\tensor5}\to\FF_2$ be the
non-bounding Barr cocycle described by Whitehouse \cite{whitehouse}.
Its composite with $\ZZ[x]^{\tensor5}\to\FF_2[x]^{\tensor5}$ is again a cocycle
and satisfies the same invariance properties, and if this composite
were a coboundary then $c$ would be too.
So $Barr^5(\ZZ[x];\FF_2)$, which is the cohomology in degree 5 of Grillet's
complex for the free commutative monoid $\NN$ with coefficients in
$\alpha^*\FF_2$, is nontrivial.

\end{document}